\documentclass[reqno,12pt]{amsart}
\usepackage{amscd,amssymb}

\theoremstyle{plain}
\newtheorem{Thm}[subsection]{Theorem}
\newtheorem{Cor}[subsection]{Corollary}
\newtheorem{Lem}[subsection]{Lemma}
\newtheorem{Prop}[subsection]{Proposition}
\newtheorem{Conj}[subsection]{Conjecture}

\theoremstyle{definition}
\newtheorem{Def}[subsection]{Definition}

\theoremstyle{remark}

\newtheorem{Rem}[subsection]{Remark}

\newtheorem{conj}{Conjecture}

\numberwithin{equation}{section}
\renewcommand{\rm}{\normalshape}

\newif\ifShowLabels
\ShowLabelstrue
\newdimen\theight
\def\TeXref#1{%
        \leavevmode\vadjust{\setbox0=\hbox{{\tt
                \quad\quad  {\small \rm #1}}}%
        \theight=\ht0
        \advance\theight by \lineskip
        \kern -\theight \vbox to
        \theight{\rightline{\rlap{\box0}}%
        \vss}%
        }}%

\ShowLabelsfalse% comment this out if labels should be printed

\renewcommand{\sec}[2]{\section{#2}\label{S:#1}%
        \ifShowLabels \TeXref{{S:#1}} \fi}
\newcommand{\ssec}[2]{\subsection{#2}\label{SS:#1}%
        \ifShowLabels \TeXref{{SS:#1}} \fi}

\newcommand{\refss}[1]{Section ~\ref{SS:#1}}

\newcommand{\reft}[1]{Theorem ~\ref{T:#1}}
\newcommand{\refl}[1]{Lemma ~\ref{L:#1}}

\newcommand{\refc}[1]{Corollary ~\ref{C:#1}}

\newcommand{\refe}[1]{\eqref{E:#1}}

\newenvironment{thm}[1]%
        { \begin{Thm} \label{T:#1}  \ifShowLabels \TeXref{T:#1} \fi }%
        { \end{Thm} }

\renewcommand{\th}[1]{\begin{thm}{#1} \sl }
\renewcommand{\eth}{\end{thm} }

\newenvironment{lemma}[1]%
        { \begin{Lem} \label{L:#1}  \ifShowLabels \TeXref{L:#1} \fi }%
        { \end{Lem} }
\newcommand{\lem}[1]{\begin{lemma}{#1} \sl}
\newcommand{\elem}{\end{lemma}}

\newenvironment{propos}[1]%
        { \begin{Prop} \label{P:#1}  \ifShowLabels \TeXref{P:#1} \fi }%
        { \end{Prop} }
\newcommand{\prop}[1]{\begin{propos}{#1}\sl }
\newcommand{\eprop}{\end{propos}}

\newenvironment{corol}[1]%
        { \begin{Cor} \label{C:#1}  \ifShowLabels \TeXref{C:#1} \fi }%
        { \end{Cor} }
\newcommand{\cor}[1]{\begin{corol}{#1} \sl }
\newcommand{\ecor}{\end{corol}}

\newenvironment{defeni}[1]%
        { \begin{Def} \label{D:#1}  \ifShowLabels \TeXref{D:#1} \fi }%
        { \end{Def} }
\newcommand{\defe}[1]{\begin{defeni}{#1} \sl }
\newcommand{\edefe}{\end{defeni}}

\newenvironment{remark}[1]%
        { \begin{Rem} \label{R:#1}  \ifShowLabels \TeXref{R:#1} \fi }%
        { \end{Rem} }
\newcommand{\rem}[1]{\begin{remark}{#1}}
\newcommand{\erem}{\end{remark}}

\newenvironment{conjec}[1]%
        { \begin{Conj} \label{Co:#1}  \ifShowLabels \TeXref{Co:#1} \fi }%
        { \end{Conj} }
\renewcommand{\conj}[1]{\begin{conjec}{#1} \sl }
\newcommand{\econj}{\end{conjec}}

\newcommand{\eq}[1]%
        { \ifShowLabels \TeXref{E:#1} \fi
           \begin{equation} \label{E:#1} }
\newcommand{\eeq}{ \end{equation} }

\newcommand{\prf}{ \begin{proof} }
\newcommand{\epr}{ \end{proof} }

\newcommand\alp{\alpha}

\newcommand\del{\delta}         \newcommand\Del{\Delta}

\newcommand\kap{\kappa}
\newcommand\lam{\lambda}

\newcommand\calD{{\mathcal{D}}}

\newcommand\calF{{\mathcal{F}}}
\newcommand\calG{{\mathcal{G}}}
\newcommand\calH{{\mathcal{H}}}

\newcommand\calL{{\mathcal{L}}}

\newcommand\calR{{\mathcal{R}}}
\newcommand\calS{{\mathcal{S}}}
\newcommand\calT{{\mathcal{T}}}

\newcommand\QQ{\mathbb{Q}}

\renewcommand\AA{\mathbb{A}}

\newcommand\GG{\mathbb{G}}

 \newcommand\grg{{\mathfrak{g}}}

 \newcommand\grp{{\mathfrak{p}}}
 \newcommand\grq{{\mathfrak{q}}}

 \newcommand\gru{{\mathfrak{u}}}

\newcommand\sdp{\times \hskip -0.3em {\raise 0.3ex
\hbox{$\scriptscriptstyle |$}}} % semidirect product

\newcommand\oB{{\overline{B}}}

\newcommand\op{{\overline{p}}}
\newcommand\oP{{\overline{P}}}
\newcommand\oq{{\overline{q}}}

\newcommand\ou{{\overline{u}}}
\newcommand\oU{{\overline{U}}}

\newcommand\tilj{{\widetilde{j}}}

\newcommand\x{\times}
\newcommand\ten{\otimes}

\newcommand\qlb{{\overline \QQ}_l}

\newcommand\Spr{{\bf Spr}}

\newcommand\Av{\text{Av}}
\newcommand\qlbt{\qlb[1](\frac{1}{2})}

\begin{document}
\title{Some results about geometric Whittaker model}
\author{Roman Bezrukavnikov, Alexander Braverman and Ivan Mirkovic}
\address{R.~B.: Department of Mathematics, University of Chicago and Clay Mathematics Institute}
\address{A.~B.: Department of Mathematics, Harvard University}
\address{I.~M.: Department of Mathematics, University of Massachusetts at
Amherst}
\thanks{All the authors were partially supported by the NSF}
\begin{abstract}
Let $G$ be an algebraic reductive group over a field of positive characteristic.
Choose a parabolic subgroup $P$ in $G$ and denote by $U$ its unipotent
radical. Let $X$ be a $G$-variety.
The purpose of this paper is to give two examples of a situation in which
the functor of averaging of $\ell$-adic sheaves on $X$ with respect to a character
$\chi:U\to \GG_a$ commutes with Verdier
duality. In the first example we take $X$ to be an arbitrary $G$-variety and
we prove the above property for all $\oP$-equivariant
sheaves on $X$ where $\oP$ is an opposite
parabolic subgroup assuming $\chi$ satisfies a strong nondegeneracy condition (such
a $\chi$ exists for some but not all choices of $P$).  In the case when $P$ is a Borel subgroup it is enough to require that the sheaf in question is 
$\oU$ equivariant where $\oU$ is the unipotent radical of $\oP$. In the second example we take $X=G$ where $G$ acts by left translations
and we prove the
corresponding result when $P$ is a Borel subgroup
for sheaves equivariant under the adjoint
action of $G$ (the latter result was conjectured by B.~C.~Ngo who proved it
for $G=GL(n)$). As an application %of the proof of the first statement
we reprove a theorem of N.~Katz and G.~Laumon about local acyclicity of
the kernel of the Fourier-Deligne transform.
\end{abstract}
\maketitle

\sec{}{Introduction}
\ssec{}{}In this paper $k$ will be an algebraically closed field of
characteristic $p>0$.
We choose a prime number $\ell$ which is
different from $p$. By a sheaf on a $k$-scheme $S$ we mean an $\ell$-adic
etale sheaf. We denote by $\calD^b(S)$ the bounded derived category of such sheaves.
For a complex $\calF\in\calD^b(S)$ we denote by $^pH^i(\calF)$
its $i$-th perverse cohomology.
Recall that for any finite subfield $k'\subset k$ and any
non-trivial character $\psi:k'\to\qlb$ we can construct the
Artin-Schreier sheaf $\calL_\psi$ on $\GG_{a,k}$.

Let $G$ denote a connected %split 
reductive group over $k$.
We shall assume that $p>h$ %is sufficiently large (with respect to $G$)
where $h$ is the Coxeter number, so that for every unipotent subgroup $U\subset G$ with Lie algebra
$\gru$ the exponential
map $\gru\to U$ is well defined and is an isomorphism. Under this condition $\grg$ carries an invariant bilinear form $(\ ,\ )$.

Let
$m:G\x G\to G$ be the multiplication map. For every
$\calF,\calG\in \calD^b(G)$ we shall denote by
$\calF\star\calG$ their "!"-convolution; in other words
\eq{shrikconvolution}
\calF\star\calG=m_!(\calF\boxtimes \calG).
\end{equation}
Similarly, we shall denote by $\calF *\calG$ the "*"-convolution of $\calF$ and
$\calG$, i.e.
\eq{starconvolution}
\calF *\calG=m_*(\calF\boxtimes \calG).
\end{equation}
%----------------------------------------------------------------------
\ssec{}{Generic characters}
%%%                         Let $U$ be a maximal unipotent subgroup of $G$. For every simple root $\alp$ of $G$ we have a one-dimensional subgroup $U_\alp\subset U$ (which is non-canonically isomorphic to $\GG_a$). We denote by $\gru_\alp\subset \gru=\text{Lie}(U)$ the corresponding 1-dimensional Lie subalgebra. We say that a homomorphism $\chi:U\to\GG_m$ is non-degenerate if the restriction of $\chi$ to every $U_\alp$ is non-trivial.

Let $P\subset G$ be a parabolic subgroup
of $G$ with a Levi decomposition
$U_P\cdot L$. Let $\oP=U_\oP \cdot L$ be an opposite parabolic.

The quotient of $U$ modulo commutant is isomorphic to 
a vector group ${\mathbb G}_a^r$  canonically up to a linear automorphism: a choice of a maximal 
torus in $L$ yields  an isomorphism $\prod\limits_{\alp\in
\Del_T^{\text{min}}(U_P)
%%%	\Del_T(U_P)^{\text{min}}
} {\mathbb G}_a\cong U/[U,U]$; since any two maximal tori in $L$ are conjugate
while the action of $L$ on ${\mathbb G}_a^r$ is easily seen to be linear, the linear structure does not depend on the maximal torus. All homomorphisms $U\to {\mathbb G}_a$
considered in the paper will be assumed to come from a {\em linear} map ${\mathbb G}_a^r \to {\mathbb G}_a$.

 We say that a homomorphism
$\chi:U_P\to \GG_a$ is non-degenerate if
for any parabolic $Q$ conjugate to $\oP$ and such that $Q\cap U_P\ne \{1\}$ the restriction $\chi|_{Q\cap U_P}$
is nontrivial.

We have a (uniquely defined) element $e_\chi\in \gru_\oP=Lie (U_\oP)$ such that $(e_\chi, x)=d\chi(x)$ for $x\in \gru=Lie(U_P)$

Recall that an element in the radical $\gru_Q$ of a parabolic subalgebra $\grq=Lie(Q)$ is called {\em Richardson} if its $Q$-orbit is open in
$\gru_Q$. It is well known  \cite[Theorem 1.3]{LS}, \cite[Theorem 7.1.1]{CM},
that in this case the neutral component of the centralizer $Z(e)^0$ is contained in $Q$ and also
that the preimage of $e$ under the moment map $\mu_Q:T^*(G/Q)\to \grg^*=\grg$ is finite. 
We say that $e\in \gru_Q$ is {\em birationally Richardson} if $Z(e)\subset  Q$. This is equivalent to the map $\mu_Q$ 
being birational onto its image and also to the condition $|\mu_Q^{-1}(e)|=1$.

\lem{} \label{L13}
The character $\chi$ is non-degenerate if and only if the element $e_\chi\subset \gru_\oP$ is birationally Richardson.
\elem

\prf
Suppose  that $e_\chi$ is birationally Richardson and let $Q$ be a parabolic conjugate to $\oP$ such that $Q\cap U_P\ne \{1\}$.
We need to show that $\chi|_{Q\cap U_P}$ is nontrivial. This property is invariant under conjugating $Q$ by an element of $U_P$.
It is well known that every $U_P$-orbit on a partial flag manifold contains a $C_L$ fixed point where $C_L$ is the center of $L$.
Thus we can assume without loss of generality that $Q\supset C_L$. Then we have $\grq = \grq\cap \gru_P \oplus \grg \cap \bar{\grp}$.
If $\chi| _{Q\cap U}$ is trivial then $e_\chi$ is orthogonal to $\grq\cap \gru$. It is also orthogonal to $\grq \cap \bar{\grp}$ by virtue of lying
in $\gru_\oP$. Thus $e_\chi$ is orthogonal to $\grq$, i.e. it is contained in its radical which contradicts the assumption that it does not
lie in any subalgebra conjugate to but different from $\gru_\oP$.

Suppose now that $e_\chi$ is not birationally Richardson. If it is not Richardson then parabolic subalgebras
conjugate to $\gru_{\oP}$ and containing $e_\chi$ correspond to points of a closed subvariety $(G/\oP)_\chi$ in $G/\oP$ of
positive dimension. Being projective, such a variety can not be contained in the affine open orbit of $U$ on $G/\oP$. 
If  $x\in (G/\oP)_\chi$ is not contained in the open orbit then the stabilizer $Q$ of $x$ satisfies $\chi|_{Q\cap U_P}=0$. 

Assume now that $e_\chi$ is Richardson but not birationally Richardson. 
Then there exists a parabolic $Q\ne P^-$ conjugate to $P^-$ and such that $e_\chi$ is in the radical of $\grq$.
It remains to show that $Q\cap U$ is nontrivial. If that was not the case then there exists $u\in U_P$ such that 
$Q=u\oP u^{-1}$. Then we claim that $Ad(u) e_\chi=e_\chi$. To check this it suffices to check that $(e_\chi,x)=(Ad(u)e_\chi,x)$
for every $x\in \grg$. If $x\in \gru_P$ this is clear since $e_\chi$ is orthogonal to the commutator of $\gru_P$. If $x\in \grq$
then both sides vanish since   $e_\chi$ is in the radical of $\grq$. Since $\grg=\gru_P\oplus \grq$ we see that $e_\chi$ is invariant
under $Ad(u)$. On the other hand, since $e_\chi$ is Richardson its stabilizer in $U$ is trivial, which is a contradiction. 
\epr

\medskip
\noindent
{\it Remark.} 
We now briefly discuss examples and non examples of the above situation.

0. If $P=B$ is a Borel subgroup and $\chi$ does not vanish on each simple root subgroup then $\chi$ is clearly nondegenerate.

1. Let $(e,f,h)$ be an $sl(2)$ triple in $\grg$ where the nilpotent $e$ is assumed to be even. Suppose that $\grp=\grg_{\leq 0}$, the sum
of eigenspaces of $h$ with nonpositive eigenvalues (we assume that the characteristic $p$
is large enough so that the $sl(2)$ triple comes from a homomorphism $SL(2)\to G$ and
all weights of $SL(2)$ in the adjoint representation are less than $p$). Then it is clear that there exists a character $\chi$
of $U_P$ such that $e=e_\chi$. It is well known that $e$ is birationally Richardson in this case.

2. The nilpotent radical of any parabolic contains Richardson elements, however,  Richardson
elements coming from a character of the radical of opposite parabolic (i.e. orthogonal to the commutant
of the radical of an opposite parabolic) exist in some but not all cases. For example,\footnote{The example
was explained to us by Zhiwei Yun.} let $G=SL(n)$
and let $P$ be the group of block upper triangular matrices with blocks of sizes $n_1, \dots, n_d$.
A Richardson element of this sort exists if and only if we have 
$n_1\leq n_2 \cdots \leq n_i\geq n_{i+1} \cdots \geq n_d$ for some $i\in [1..d]$. 
To see this notice that an element $e$ is Richardson if and only if the Young diagram $\lambda$
corresponding to its Jordan type is transposed to the diagram  $\mu$ of the partition $n=n_1+\cdots + n_d$.
In particular, this implies that $\dim (Ker(e)) = \max \{ n_i\}$. However, $\dim (Ker(e)) \geq n_d+ \sum\limits _{i, \, n_i> n_{i+1}} n_i-n_{i+1}$
which is greater than $\max \{ n_i\}$ unless $n_i$ satisfy the above inequalities. This shows the inequalities are necessary
for existence of a Richardson element of the specified form. To see they are sufficient consider the homogeneous
ideal  $J$ in $k[x,y]$ corresponding to the partition $m_1=\cdots = m_i=n_i$, $m_j=n_j$ for $j>i$; thus monomials
$x^iy^j$, $j\leq m_i$ form a basis in $k[x,y]/J$. The subspace $V\subset k[x,y]/J$ spanned by monomials $x^iy^j$, $j\leq n_i$
is clearly invariant under multiplication by $x$. The operator of multiplication by $x$ is easily seen to have the Jordan form
given by the diagram $\lambda$,  
 while decomposing $V=\oplus V_s$ where $V_s$ is spanned by monomials $x^s y^j$ we see that multiplication by $x$ corresponds to a matrix
 with nonzero blocks right above block diagonal as required.

Recall the well known fact  that if $G$ is of type $A$ then every Richardson element is birationally Richardson,
thus the previous paragraph provides examples of non-degenerate characters of $U_P\subset SL(n)$.

3. Let $G=Sp(6)$ and let $P$ be the stabilizer of a line in the tautological representation. 
Then $U_P$ is the 5-dimensional Heisenberg group. For any character $\chi$ not vanishing
on the 1-dimensional center of $U_P$ the element $e_\chi$ is Richardson but not birationally
Richardson.

%-----------------------------------------------------------------------------
\ssec{dickens}{}
Let $X$ be a $G$-variety. Assume that  $U$ is a  subgroup of $G$ and
$\chi:U\to\GG_a$ is a homomorphism. Let $a:U\x X\to X$ denote the action map and
let $p:U\x X\to X$ be the projection to the second multiple. Let
$\Av_{U,\chi,*}:\calD^b(X)\to\calD^b(X)$ be the functor
sending  $\calF\in\calD^b(X)$ to
$a_* (\chi^*\calL_\psi\boxtimes\calF) \left( \qlbt \right) ^{\ten \dim U}$.
Similarly we define
the functor $\Av_{U,\chi,!}$ by replacing $a_*$ by $a_!$.
We have the natural morphism $\Av_{U,\chi,!}\to \Av_{U,\chi,*}$.
The main result of this paper is the following:
%-----------------------------------------------------------------------------------
\th{dickens}
Let $U\subset G$ be the unipotent radical of a parabolic subgroup
$P\subset G$ and let $\chi:U\to\GG_a$ be non-degenerate.
\begin{enumerate}
\item
Let $\oP$ denote a parabolic subgroup of $G$ opposite to $P$ and let $\oU$ denote
its unipotent radical.
Let $\calF\in\calD^b(X)$ be $\oP$-equivariant. Then
the natural morphism
$$
\Av_{U,\chi,!}\calF\to \Av_{U,\chi,*}\calF
$$
is an isomorphism.
\item
Assume that $P$ is a Borel subgroup. Then the above assertion is true for all $\oU$-equivariant sheaves $\calF$ (not just $\oP$-equivariant ones).
\item
Let $\calF\in\calD^b(G)$ be equivariant with respect to the adjoint action and assume that $P$ is a Borel subgroup.
Then the natural morphism
$$
\Av_{U,\chi,!}\calF\to \Av_{U,\chi,*}\calF
$$
is an isomorphism.
\end{enumerate}
\eth

\noindent
{\it Remarks.}

\noindent
0.
In the above cases the averaging functors preserve  perversity:\
the functor is defined as a direct image under an affine morphism, thus if $\calF$ is perverse then
$\Av_{U,\chi,!}\calF$ is in perverse degrees $\ge 0$, while
 $\Av_{U,\chi,*}\calF$ is in perverse degrees $\le 0$. Since the two complexes are isomorphic,
 they are perverse.

\noindent
1. The second statement of \reft{dickens} was communicated to the
first author as a conjecture by B.~C.~Ngo who also proved it for $G=GL(n)$.

\noindent
2. In the case $G=GL(n)$ a (much more involved) analogue of \reft{dickens}(3)
is used in \cite{Ga} (Theorem 5.1) in order to complete the proof of the geometric
Langlands conjecture for $GL(n)$. 
%We believe that both statements of \reft{dickens}
%might have something to do with a possible generalization of Theorem 5.1 of
%\cite{Ga} to the case of arbitrary reductive group.

\noindent
3. In the next section we also explain how the main step in the proof of
\reft{dickens}(1) allows to reprove one of the main results of \cite{KL}.

\noindent
4. \reft{dickens} also holds, with a parallel proof, when $k$ is an algebraically closed field of
characteristic 0 and $\ell$-adic sheaves are replaced by holonomic
$\calD$-modules (in this case one has to replace $\calL_\psi$ by the
$\calD$-module corresponding to the function $e^x$).

\noindent
5. In the previous version of this paper \cite{publ} the statement of \reft{dickens}(1) was formulated  for an arbitrary parabolic subgroup $P$ of $G$  under 
the weaker assumption that $\calF$ is $\oU$-equivariant and with a different definition of a generic character. However, our proof of that assertion was wrong (unless $P$ is a Borel subgroup) and the definition of a generic character irrelevant 
as was recently pointed out to us by Z.~Yun. In fact, the statement itself is wrong for any choice of the character as the following counterexample shows. 
Suppose that $H\subset G$ is a subgroup such that
\begin{enumerate}
\item $\dim(H\cap U)>0$.
\item $H\cap U\subset Ker(\chi)$.
\item $\oU\subset H$.
\end{enumerate}
If $X=G$ and $\calF$ is the constant sheaf on $H\subset G$ then 
the stalks of  $\Av_{U,\chi,*}(\calF)$ and $\Av_{U,\chi,!}(\calF)$ at $1\in G$ are isomorphic
to $H^* \left( H\cap U \right) [ \dim U] \left( \frac{\dim U}{2} \right)$ and 

\noindent
$H^*_c\left(H\cap U \right)   [ \dim U] \left( \frac{\dim U}{2}\right)$
respectively, thus the two averaging functors produce nonisomorphic complexes.
An example of such a subgroup is as follows. Consider $G=SL(3)$ and let $P$ and $\oP$ be the groups of block upper triangular 
and block lower triangular matrices respectively with blocks of sizes $(2,\, 1)$. Let $H$ be the group of block lower triangular matrices
with blocks of sizes $(1,\, 2)$. Then conditions (1) and (3) are clear. Also notice that in this case $U$ is abelian and all nontrivial characters 
of $U$ are conjugate under the action of $L=P/U$, thus for any choice of the character (2) also holds
possibly after replacing $H$ by its conjugate by an element of $L$.

\noindent
6. On the other hand, the present formulation of  \reft{dickens}(1) may not be optimal:
it was pointed out to us by Z.~Yun that the statement may hold under a weaker assumption of the character. However, proving 
the more general statement requires new ideas and is beyond the scope of this document.
  
\noindent  
7. Originally \reft{dickens}(3) was formulated for any $P$, it was deduced from the above (stronger and incorrect) formulation of \reft{dickens}(1). 

We do not know if the latter assertion is true. On the other hand, we believe that the following generalization of \reft{dickens}(3) holds. Assume for simplicity that $p$ is sufficiently large, let $U$ be a (unipotent) Premet group associated to a nilpotent orbit in $\grg^*$ (cf. e.g. \cite{Premet}). It comes with a natural character $\chi:U\to \GG_m$.
\conj{premet}
\begin{enumerate}
\item
Then the assertion of \reft{dickens}(3) holds for any Premet pair $(U,\chi)$.
\item
Let $U$ and $\chi$ be
as above.
Then for any  irreducible perverse sheaf $\calF\in\calD^b(G)$
equivariant with respect to the adjoint action,
$\Av_{U,\chi,!}\calF$ is an irreducible perverse sheaf
or zero.
\end{enumerate}
\econj

%----------------------------------------------------------------------------------------------
\sec{}{Proof of \reft{dickens}(1)}
\ssec{}{Cleanness}
Let $Z$ be an algebraic variety over $k$ and let $j:Z_0\to Z$ be an open
embedding. We shall say that $\calG\in\calD^b(Z_0)$ is {\it clean with respect to
$j$} if the natural map $j_!\calG\to j_*\calG$ is an isomorphism.

Let $X$ be any $\oP$-variety. Consider the variety
$G\underset{\oP}\x X$. We have the natural open embedding
$j:U\x X\to G\underset{\oP}\x X$. We will prove
\reft{dickens}(1) by a series of reductions.
We claim that \reft{dickens}(1)  follows from
%--------------------------------------------------------------------------------------
\th{cleanness}
Let $X$ be a $\oP$-variety and let
$\calF\in\calD^b(X)$ be $\oP$-equivariant. Then the sheaf
$\chi^*\calL_\psi\boxtimes \calF$ is clean with respect to $j$.
In other words the natural morphism
\eq{}
j_!(\chi^*\calL_\psi\boxtimes \calF)\to j_*(\chi^*\calL_\psi\boxtimes\calF)
\end{equation}
is an isomorphism. The same thing is true under the assumption that $\calF$ is $\oU$-equivariant if $P$ is a Borel subgroup.
\eth

%-------------------------------------------------------------------------------------
\ssec{}{
\reft{cleanness}
implies
\reft{dickens}(1)
}
Indeed if $X$ is a $G$-variety then we have the natural proper map $b:G\underset{\oP}\x X\to X$
sending every $(g,x)\, \text{mod}\, \oP$ to $g(x)$. Moreover, we have $b\circ j=a$ (recall that
$a:U\x X\to X$ denotes the action map). Hence \reft{cleanness} and the fact that $b$ is
proper imply that
$$
\Av_{U,\chi,!}\calF=b_!(j_!(\chi^*\calL_\psi\boxtimes \calF))=
b_*(j_*(\chi^*\calL_\psi\boxtimes \calF))=\Av_{U,\chi,*}\calF.
$$
It remains to prove \reft{cleanness}. Note that in the formulation of \reft{cleanness} we do not need
$X$ to be a $G$-variety but only a $\oP$-variety.
%------------------------------------------------------------------------------------

\ssec{}{
A reformulation  of the
\reft{cleanness}
}
Let $\pi:G\x X\to G\underset{\oP}\x X$ be the natural projection. Also let
$\tilj:U\cdot\oP\x X\to G\x X$ be the natural embedding. It follows
from the smooth base change theorem  that it is enough to show that the
natural map
$$
\tilj_!\pi^*(\chi^*\calL_\psi\boxtimes \calF)\to \tilj_*\pi^*(\chi^*\calL_\psi\boxtimes\calF)
$$
is an isomorphism (note that we have the natural identification $U\cdot\oP\x X$ with
$\pi^{-1}(U\x X)$).

The sheaf $\pi^*(\chi^*\calL_\psi\boxtimes \calF)$ is obviously
$(U,\chi)$-equivariant with respect to the $U$-action by multiplication on
the left. We claim that it is also $\oP$-equivariant with respect to
multiplication on the right, i.e. with respect to the $\oP$-action
on $U\cdot \oP\x X$ given by $\oq:(u,\op,x)\mapsto (u,\op\oq,x)$.
Indeed, the map $\pi$ from $U\cdot\oP\x X$ to $U\x X$ is given
by $\pi:(u,\op,x)\mapsto (u, \op(x))$ (since the action of
$\oP$ on $G\x X$ is given by $\op:(g,x)\mapsto (g\op^{-1},\op x)$).
Thus
$$
\pi(u,\op\oq,x)=(u,\op\oq(x))=(u, \op\oq\op^{-1}(\op(x)))
$$
and our statement follows from $\oP$-equivariance of $\calF$.

Similar argument shows that if $\calF$ is only $\oU$-equivariant, then $\pi^*(\chi^*\calL_\psi\boxtimes \calF)$ is $\oU$-equivariant with respect to the right multiplication.

Hence we see that \reft{cleanness} follows from the following lemma.
\lem{}
Consider the action of $U\times \oP
$
on $U\cdot \oP\subset G$
given by left and right multiplications.
For any variety $X$, if
$\calG\in\calD^b(U\cdot\oP\x X)$ is $(U,\chi)$-equivariant on the left
and $\oP$-equivariant on the right, then
the natural map given by the inclusion
$\tilj:U\cdot\oP\x X\to G\x X$,
$$
\tilj_!\calG\to\tilj_*\calG
,$$
is an
isomorphism. Similar assertion holds if $\calF$ is only $U\times \oU$ equivariant if $P$ is a Borel subgroup.
\elem
\prf
Let $Z$ denote the complement of $U\cdot\oP$ in $G$
and let $i$ be the natural embedding of $Z\x X$ to $G\x X$.
Since $\tilj_*\calG$ is also $(U,\chi)$-equivariant on the left
and $\oU$-equivariant on the right it is enough to show that
for every complex $\calH$ on $G\x X$ with the above equivariance properties
we have $i^*\calH=0$. However, it is clear that this follows from:
%%% the following
%%%\epr
%-------------------------------------
\lem{stabilizer}
\begin{enumerate}
\item
Let $g\in Z$. Let $S_g\subset U\x \oP$ denote the set of all
pairs $(u,\op)$ such that $ug\op=g$. Let also $U_g$ be the projection of $S_g$
to $U$. Then the restriction of $\chi$ to $U_g$ is non-trivial.
\item
Similarly, assume that $P$ is a Borel subgroup. Let $g\in Z$. Let $S'_g\in U\x \oU$ denote the set of all
pairs $(u,\ou)$ such that $ug\ou=g$. Let also $U_g'$ be the projection of $S_g'$
to $U$. Then the restriction of $\chi$ to $U_g'$ is non-trivial.
\end{enumerate}
\elem
\prf
The first assertion is clear from the definition of a generic character.

%The 2nd assertion is standard (and its proof is essentially a word-by-word repetition of the above argument).
If $P=B$ then Theorem \ref{dcikens}(1) applies, see Remark 0 after Lemma \ref{13}.
Also, in this case $\oU$ coincides with the set of unipotent elements in $\oP$, thus we have: $S_g=S_g' $ and $U_g=U_g'$. Thus (2) follows from (1).
%
%\footnote{The difference between the case when $P$ is a Borel subgroup and the case of general $P$ is that in the former case we always
%have $U\cap \tilw\oP\tilw^{-1}=U\cap \tilw\oU\tilw^{-1}$, and the same formula is wrong in the general case.}
\epr
\epr
%%%\epr
%--------------------------------------------------------------------------------
\cor{cleanP}
Let $U$ be the unipotent radical of a Borel subgroup $B$.
Let $j$ denote the open embedding of $B$ into $G/\oU$.
Let $\calF$ be any $(U,\chi)$-equivariant sheaf on $B$
(with respect to the left multiplication action).
Then the natural morphism
$$
j_!\calF\to j_*\calF
$$
is an isomorphism. In other words, every $(U,\chi)$-equivariant sheaf on
$B$ is clean with respect to $j$.
\ecor
%------------------------------------------------------------------
\prf
Let $L$ be the Levi factor of $B$ (since $B$ is now assumed to be a Borel subgroup, $L$ is a maximal torus of $G$). The isomorphism $L\simeq \oB/\oU$ gives rise
to a natural action of $\oB$ on $L$. Since the action of $\oU$ on $L$ is trivial
it follows that every $\calF\in\calD^b(L)$ is automatically
$\oU$-equivariant.

We have the natural identifications $U\x L\simeq B$ (by multiplication map) and
 $G\underset{\oB}\x L\simeq G/\oU$ (sending every $(g,l)\mod \oB$ to $gl\mod \oU$). Under
these identification the embedding $j:B\to G/\oU$ becomes equal to the
natural embedding $U\x L\to G\underset{\oB}\x L$ considered in
\reft{cleanness} (for $X=L$). Also the fact that $\calF$ is
$(U,\chi)$-equivariant implies that as a sheaf on $U\x L$ it can be decomposed as
$\calF=\chi^*\calL_\psi\boxtimes\calF'$ for some $\calF'\in\calD^b(L)$.
Hence \refc{cleanP} is a particular case of \reft{cleanness}.
\epr
%---------------------------------------------------------------------------
\ssec{}{Application to Katz-Laumon theorem}
Consider the variety $\AA^1\x \GG_m$ with coordinates
$(x,y)$. Let $f:\AA^1\x\GG_m\to \AA^1$ be given by $f(x,y)=\frac{x}{y}$. Let
also $i:\AA^1\x\GG_m\to \AA^2$ denote the natural embedding and let $\pi:\AA^1\x \GG_m\to \GG_m$
be the projection to the second variable. The following theorem is proved in
\cite{KL}.
%----------------------------------------------------------------
\th{katz-laumon}
For every $\calF\in\calD(\GG_m)$ the natural map
\eq{kl}
i_!(f^*\calL_\psi\ten\pi^*\calF)\to i_*(f^*\calL_\psi\ten\pi^*\calF)
\end{equation}
is an isomorphism.
\eth
%---------------------------------------------------------------
Below we explain that \reft{katz-laumon} may be viewed as a particular case
of \refc{cleanP}.
\prf
Take now $G=SL(2)$ and let $B$ and $\oB$ be respectively the subgroups of
lower-triangular and
upper-triangular matrices, with unipotent radicals $U$ and $\oU$.
We denote the natural isomorphism between $U$ and $\GG_a$ by $\chi$.

Let us identify $G/\oU$ with $\AA^2\backslash\{0\}$ by
$g\oU\mapsto
g(e_1)$ for the first standard basis vector $e_1$ of $\AA^2$.
Then
$B\subset G/\oU$ is identified with
$\AA^1\x \GG_m\subset\ \AA^2\backslash\{0\}$ by
$
\begin{pmatrix}
\lam& 0\\
t &\lam^{-1}
\end{%%%small
pmatrix}
\leftrightarrow
(t,\lam^{-1})
$.
The sheaf
$f^*\calL_\psi\ten\pi^*\calF$ is
$(U,\chi)$-equivariant, so by
\refc{cleanP} this sheaf is clean for the embedding
$\AA^1\x \GG_m\subset\ \AA^2\backslash\{0\}$.
It remains to observe that the resulting sheaf
on
$\AA^2\backslash\{0\}$ is clean for the embedding into
$\AA^2$ since the cone of the canonical map between the shriek and
star direct images is zero -- it is a
$(U,\chi)$-equivariant sheaf supported at  a point  $\{0\}$.
\epr
%-----------------------------------------------------------------------------------------------------------
%%%%%%%%%%%%%%%%%%%%%%%%%%%%%%%%%%%%%%%%%%%%%%%%%%%%%%%%%%%%%%%%%%%%%%%%%%%%%%%%%%%%%%%%%%%%%%%%%%%%%%%%%%%%%
%%%%%%%%%%%%%%%%%%%%%%%%%%%%%%%%%%%%%%%%%%%%%%%%%%%%%%%%%%%%%%%%%%%%%%%%%%%%%%%%%%%%%%%%%%%%%%%%%%%%%%%%%%%%%
%-------------------------------------------------------------------------------------------------------------

\sec{}{Proof of \reft{dickens}(3)}
In this Section we will eventually need to assume that $P$ is a Borel subgroup. However, let us start working with arbitrary $P$ and make the above assumption later.
%----------------------------------------------------
\ssec{}{Horocycle transform}Let
$P$ be a parabolic subgroup in $G$ and let $Y_P$ denote the
variety of all parabolic subgroups of $G$ which are conjugate to
$P$. We also denote by $W_P$ the variety of $P$-horocycles,
i.e., the pairs $(Q\in Y_P,x\in
G/U_{Q})$ where $U_{Q}$ denotes the unipotent radical of $Q$
(see section \refss{WQ} below for a more direct definition of $W_P$).
We have the natural map $p:W_P\to Y_P$.

%%%%%%%%%%%%%%%%%%%%%%% do sih ispravil, chort vozmi the rest

We also have the natural morphisms
$\alp:G\x Y_P\to G$ and
$\beta:G\x Y_P\to W_P$ where $\alp$ is just the projection to the
first multiple and $\beta$ sends
%%%every
$(g,Q)$ to
$(Q,g \mod U_{Q})$. We define two functors
$\calR_P:
\calD^b(G)\to\calD^b(W_P)$ and $\calS_P:\calD^b(W_P)\to\calD^b(G)$
by setting
$$
\calR_P(\calF)=\beta_!\alp^*(\calF)\ten\left( \qlbt \right)^{\ten \dim U_P}
$$
and
$$
\calS_P(\calG)=\alp_!\beta^*(\calG)\ten\left( \qlbt \right)^{\ten \dim U_P}.
$$
The following lemma is proved in \cite{MV} when $Q$ is  a Borel subgroup in
$G$.
\lem{directsummand}
The identity functor is a direct summand of $\calS_P\circ\calR_P$.
\elem
\prf
Let $\calT_Q=\{(Q\in Y_P,u\in U_{Q})\}$. We have the natural map
$p_P:\calT_Q\to G$ sending every $(Q,u)$ to $u$ (clearly the image of
$p_P$ lies in the set of unipotent elements in $G$). Let
$\Spr_P=(p_P)_!\qlb[2\dim Y_P](\dim Y_P)$. It is known (cf. \cite{BM})
that $\Spr_P$ is perverse and that it contains the skyscraper sheaf
$\delta_e$ at the unit element $e\in G$ as a direct summand.
We set $\Spr_P=\del_e\oplus\Spr'_P$.

On the other hand, arguing as in \cite{MV} we can show that
for every $\calF\in\calD^b(G)$ we have a canonical isomorphism
\eq{}
\calS_P\circ\calR_P(\calF)=\calF\star\Spr_P.
\end{equation}
Hence
\eq{}
\calS_P\circ\calR_P(\calF)=\calF\oplus(\calF\star\Spr'_P)
\end{equation}
which finishes the proof.
\epr
%-----------------------------------------------------------------------------------------------
\ssec{WQ}{Another definition of $W_P$}
One can identify $W_P$ with $(G/U_P\x G/U_P)/M$ where $M=P/U$ acts on
$G/U_P\x G/U_P$ diagonally. The identification is given by the
map
$$
(x_1 \mod U_P, x_2\mod U_P)\mapsto (x_2 P x_2^{-1}, x_1x_2^{-1}
\mod U_{P})
.$$
Under this identification the natural left and right $G$-actions on
$(G/U_P\x G/U_P)/M$
give two actions of $G$ on $W_P$, which we still call the
 "left" and "right" action. The left action is just the natural $G$-action
in the fibers of $p$.
The right action is given by
$$
g:(Q,x)\mapsto (gQg^{-1}, xg^{-1}\mod gU_{Q}g^{-1}).
$$
The corresponding adjoint action is given by
$$
g:(Q,x)\mapsto (gQg^{-1}, g xg^{-1}\mod gU_{Q}g^{-1}).
$$

We now claim the following

Assume now that we could prove that for any $\calG\in\calD^b(W_{\oP})$ which is equivariant with respect to the adjoint
action and for every non-degenerate character $\chi:U\to \GG_a$ the natural map
\eq{av-p}
\Av_{U,\chi,!}\calG\to\Av_{U,\chi,*}\calG
\end{equation}
is an isomorphism. Then the same would be true for any 
 $\calF\in\calD^b(G)$  which is equivariant with respect to the adjoint action. the map
\eq{mainthing}
\Av_{U,\chi,!}\calF\to\Av_{U,\chi,*}\calF
\end{equation}
would be an isomorphism. Indeed, since by \refl{directsummand}
$\calF$ is a direct summand of $\calS_\oP\circ\calR_\oP(\calF)$
it is enough to show that \refe{mainthing} holds for the latter.
It follows from the fact that $\alp$ is a proper morphism that
we have the natural isomorphisms of functors
$$
\Av_{U,\chi,!}\circ\calS_\oP\simeq \calS_\oP\circ\Av_{U,\chi,!}\ \text{and}\
\Av_{U,\chi,*}\circ\calS_\oP\simeq \calS_\oP\circ\Av_{U,\chi,*}.
$$
Hence it is enough to show that \refe{mainthing} holds for
$\calR_\oP(\calF)$. However, it is clear that $\calR_\oP$ maps
ad-equivariant complexes to ad-equivariant ones which finishes
the proof by \refe{av-p}.

We do not know if \refe{av-p} is indeed true for arbitrary $P$. However, we claim the following:
\th{upsihorocycle}
Let $P$ be a Borel subgroup in $G$ and let $U$ be its unipotent radical.
%%%Take $Q=\oP$.
Let $\calG\in\calD^b(W_{\oP})$ be equivariant with respect to the adjoint
action.
Then for every non-degenerate character $\chi:U\to \GG_a$ the natural map
$$
\Av_{U,\chi,!}\calG\to\Av_{U,\chi,*}\calG
$$
is an isomorphism (here averaging is performed with respect to the
left action).
\eth

\ssec{}{}
The rest of this section is occupied by the proof of \reft{upsihorocycle}.

Let $Y_\oP^0$ denote the open $U$-orbit on $Y_\oP$ and let $W_\oP^0$ denote its
preimage in $W_\oP$.

First of all we claim that both $\Av_{U,\chi,*}\calG$ and
$\Av_{U,\chi,!}\calG$ are equal to the extension by
zero of their restriction to $W_\oP^0$. Indeed we must show that the
$*$-restriction of either of these sheaves to the fiber
of $p:W_\oP\to Y_\oP$ over any parabolic $Q$ which is not
opposite to $P$ is equal to zero. Let us denote this restriction
by $\calH$. This is a complex of sheaves on $p^{-1}(Q)=G/U_{Q}$.
The fact that $\calG$ is equivariant with respect to the adjoint action
implies that both $\Av_{U,\chi,*}\calG$ and
$\Av_{U,\chi,!}\calG$ are equivariant with respect to the adjoint action of
$U$. Hence $\calH$ is equivariant with respect to the left action of
$U\cap U_{Q}$. On the other hand, it is clear that $\calH$ is
$(U,\chi)$-equivariant with respect to the left action
of $U$. Thus our statement follows from the following result which is
equivalent to \refl{stabilizer}:
let $Q$ be as above (i.e. $Q$ is conjugate to $\oP$ but it is not
in the generic position with respect to $P$);
then the restriction of $\chi$ to $U\cap U_{Q}$ is non-trivial.

It remains to show that
the map $\Av_{U,\chi,!}\calG\to \Av_{U,\chi,*}\calG$ is an isomorphism
when restricted to $W_\oP^0$.

The map $u\mapsto u \oP u^{-1}$ is an isomorphism between $U$ and $Y_\oP^0$.
Let $\kap:W_\oP^0\to U$ be the composition of the natural projection
$W_\oP^0\to Y_\oP^0$ with this isomorphism. Define now a new
$G$-action on $W_\oP^0$ (denoted by $(g,w)\mapsto g\x w$) by
$$
g\x w=\kap(w)g\kap(w)^{-1}(w)
$$
(in the right hand side we use the standard left action of $G$ on $W_\oP$).

To finish the argument we need the following general (and basically
tautological) result:
%----------------------------------------------------------------------------
\lem{twoactions}
a) Let $H$ be an algebraic group, and $X$ be an algebraic variety equipped
with two actions $\phi_1$, $\phi_2$ of $H$. Suppose that the two actions
differ by a conjugation, i.e. there
exists a morphism of algebraic varieties $c: X\to H$, such that
$$\phi_1(g)(x)=\phi_2(c(x)\cdot g\cdot c(x)^{-1})(x)$$
for all $g\in H$, $x\in X$. Then for any character $\chi:H\to \GG_a$ we have
canonical isomorphisms of the averaging functors corresponding to the two
actions:
$$
Av_{H,\chi,!}^{\phi_1}=Av_{H,\chi,!}^{\phi_2},
$$
$$
Av_{H,\chi,*}^{\phi_1}=Av_{H,\chi,*}^{\phi_2}.
$$

b) Let $H_1$, $H_2$ be two algebraic groups, $\phi_i$ be an action
of $H_i$ on an algebraic variety $X_i$ (where $i=1,2$). Let $f:X_1\to X_2$
be a morphism, and
assume that
there exists a morphism $s: H_1\times X_1 \to H_2$, such that $$
f(\phi_1(h_1)(x_1))=
\phi_2(s(h_1,x_1))(f(x))$$
for $x_1\in X_1$, $h_1\in H_1$.
Then for
any $H_2$-equivariant complex of constructible
sheaves on $X_2$ the complex $f^*(X)$ is also $H_1$ equivariant.
$\square$
\elem

Part (a) of \refl{twoactions} shows
 that both averaging functors $\Av_{U,\chi,!}$ and $\Av_{U,\chi,*}$
 do not change
when we replace the old action by the new one. Also, since
 our $\calG$
is equivariant with respect to the
adjoint action it follows  that $\calG|_{W_\oP^0}$ is also
 equivariant with respect
to the new action of $\oU$ by part (b) of \refl{twoactions}.
 The statement now follows from
\reft{dickens}(2).

{\bf Acknowledgement.} We thank Zhiwei Yun for pointing out the mistake in the published version of the paper and for helpful discussions in the course of our work on its partial correction.

\end{document}